\documentclass{amsart}

\usepackage{amssymb}
\usepackage[all]{xy}
\usepackage{hyperref}
\newtheorem{thm}{Theorem}%[section]

\newtheorem{cor}[thm]{Corollary}

\newtheorem{prop}[thm]{Proposition}

\newtheorem{conj}[thm]{Conjecture}
\theoremstyle{definition}

\newtheorem{rem}[thm]{Remark}          

\newtheorem*{ack}{Acknowledgments}      % \renewcommand{\theack}{} 

\newtheorem{defn-thm}[thm]{Definition--Theorem}  %!!!!!!!!!!!!!!!!!!!!!!!!
\newtheorem{defn-lem}[thm]{Definition--Lemma}  %!!!!!!!!!!!!!!!!!!!!!!!!
\theoremstyle{remark}

\renewcommand{\a}[0]{{\mathbb A}}

\newcommand{\p}[0]{{\mathbb P}}

\newcommand{\map}[0]{\dasharrow}

\def\into{\DOTSB\lhook\joinrel\to}

\def\loccoh#1.#2.#3.#4.{H^{#1}_{#2}(#3,#4)}

\DeclareMathAlphabet{\mathchanc}{OT1}{pzc}%
                                {m}{it}

\usepackage[all]{xy}\xyoption{dvips}

\begin{document}
\bibliographystyle{amsalpha}

\hfill\today

 \title[Comment]{Comment on: Pseudo-effectivity of the relative canonical divisor ... by Zsolt Patakfalvi}
\author{J\'anos Koll\'ar}

\maketitle

We discuss a more elementary approach to the existence of 
non-uniruled ramified covers of varieties, which was one of the steps in the proof of the main theorem  in \cite{pat-psef}. This approach gives reasonably good control of the degree of the ramification divisor. A disadvantage is that we get non-uniruledness only for very general members of a family of covers, whereas 
both \cite{pat-psef} and \cite{pat-zda} give an open subset of non-uniruled members.

We work over an algebraically closed field $k$ of arbitrary characteristic.
The main observation is the following.

\begin{prop}\label{non.ur.prop}  Let $X^n\subset \p^N$ be an integral $k$-variety
and $H\subset \p^N$ a very general hypersurface of degree $\geq n+1$. 
Then $X\cap H$ is not uniruled.
\end{prop}

Proof. As a direct consequence of a result of Matsusaka \cite{mats-68}, it is sufficient to find one hypersurface
$H_0$ of degree $n+1$ such that $X\not\subset H_0$ and at least one of the irreducible components of $X\cap H_0$ is not uniruled; see \cite[IV.1.7--8]{rc-book} for details.

Choose now a linear projection  $\pi:\p^N\map \p^n$ that is  finite on $X$. By \cite[Thm.1.1]{rei-woo} there is a hypersurface $W\subset \p^n$ of degree $n+1$ that is not uniruled. We can then take $H_0$ to be the closure of $\pi^{-1}(W)$. \qed

\begin{cor} Let $Y$ be a smooth, projective $k$-variety of dimension $n$.
Then there is a non-uniruled, smooth, projective $k$-variety $Y'$ of dimension $n$, and a finite morphism $p:Y'\to Y$ of degree $n+2$. 
\end{cor}

Proof.  Choose an embedding $Y\into \p^N$ and let 
 $X\subset \p^{N+1}$ be the cone over $Y$.  Set $Y':=X\cap H$, where 
$H$ is a very general hypersurface of degree $n+2$. Then $Y'$ is  
non-uniruled by Proposition~\ref{non.ur.prop} and we can take $p$ to be the projection from the vertex of the cone.\qed
\medskip

Note that Proposition~\ref{non.ur.prop} is sharp in case $X$ is a linear subspace. However one expects to do better in all other cases. 
Mainly to encourage further study of this question, let me state the strongest variant that could be true.

\begin{conj} Let $X$ be a smooth projective variety and $|H|$ a very ample linear system on $X$ such that $K_X+H$ is effective. Then a general $H\in |H|$ is not uniruled.
\end{conj}

This is true in characteristic 0, but there is very little evidence for this in positive characteristic. It is not even clear that it holds for $\dim X=3$. Note that the smoothness of $X$ is essential here; cones over uniruled varieties with ample canonical class give singular counterexamples.

\begin{rem} The proof of \cite{rei-woo} is quite subtle, but it leads to  simple,  explicit examples defined over finite fields; see  \cite[3.11]{rei-woo}.

Still it may be worthwhile to find an elementary argument, giving possibly worse degree bounds. As the most naive example, let $\bigl(g(x, y)=0\bigr)\subset \a^2$ be a nonrational cubic. Then
$$
\bigl( g(x_1, y_1)=\cdots=  g(x_n, y_n)=0\bigr)\subset \a^{2n}
$$
is a complete intersection of degree $3^n$ that contains no rational curves.
A birational projection of it to $\p^{n+1}$ gives a non-uniruled hypersurface
of degree $3^n$.
\end{rem}

\begin{ack} I thank T.~Murayama and Zs.~Patakfalvi 
for  helpful comments and  discussions.
Partial  financial support    was provided  by  the NSF under grant number
DMS-1901855.
\end{ack}

\bigskip

  Princeton University, Princeton NJ 08544-1000, \

\email{kollar@math.princeton.edu}

\end{document}